 \documentclass[letterpaper, 10 pt , conference]{ieeeconf}
\usepackage{graphicx} % for pdf, bitmapped graphics files
\usepackage{epsfig}
\usepackage{epstopdf}
\usepackage{stmaryrd}
\usepackage{cite} % allows to put citations like [1]-[5]
\usepackage{amsmath} % assumes amsmath package installed
\usepackage{amssymb}  % assumes amsmath package installed
\usepackage{mathrsfs}
\usepackage{color}
\usepackage{nicefrac}
\usepackage[linesnumbered, ruled, vlined]{algorithm2e}
\usepackage{mathtools}
\mathtoolsset{showonlyrefs}
\usepackage{amsthm}
\graphicspath{{.}}

\theoremstyle{plain}
\newtheorem{thm}{Theorem}
\usepackage{etoolbox}
\usepackage{float}
\usepackage{tabularx}
  
\def\th@definition{
	\thm@headfont{\itshape} % Heading font is italic
	\thm@notefont{} % Note is same as heading
}
\makeatother

\usepackage{color}

\theoremstyle{definition}
\newtheorem{dfn}{Definition}

\newtheorem{lemma}{Lemma}

\newtheorem{assumption}{Assumption}

\newcommand{\R}{\mathbb{R}}

\makeatletter
\def\endthebibliography{%
	\def\@noitemerr{\@latex@warning{Empty `thebibliography' environment}}%
	\endlist
}
\makeatother

\makeatletter
\patchcmd{\@makecaption}
{\scshape}
{}
{}
{}
\makeatother

\title{\LARGE \bf
Asynchronous Parallel Nonconvex Optimization \\ Under the Polyak-{\L}ojasiewicz Condition
}

\author{Kasra Yazdani and Matthew Hale$^{\ast}$% <-this % stops a space
\thanks{$^{\ast}$Authors are with the Department of  Mechanical and Aerospace Engineering at the University of Florida, Gainesville, FL USA. Emails: \texttt{\{kasra.yazdani,matthewhale\}@ufl.edu}. This work was supported in part by a Task Order contract with the Air Force Research Laboratory, Munitions Directorate, at Eglin AFB, by AFOSR under Grant FA9550-19-1-0169, and
 by ONR under Grant N00014-19-1-2543.}
}

\begin{document}

\maketitle
\thispagestyle{empty} % Removes the page number in the first page

\begin{abstract}

Communication delays and synchronization are major bottlenecks for parallel computing, and 
tolerating asynchrony is therefore crucial for accelerating parallel computation. 
Motivated by optimization problems that do not satisfy convexity assumptions, 
we present an asynchronous block coordinate descent algorithm for nonconvex optimization problems whose objective functions satisfy the Polyak-{\L}ojasiewicz condition. 
This condition is a generalization of strong convexity to nonconvex problems and requires
neither convexity nor uniqueness of minimizers. 
%Contrary to local gradient descent
%that runs gradient descent in parallel on different processors and averages their updates over time, 
%in our algorithm, at each step, a processor updates a subset of the decision variables
%based on possibly outdated information it stores onboard,
%and then processors share these updates through asynchronous communication. 
Under only assumptions of mild smoothness of objective functions and bounded delays,
we 
% provide rules for selecting stepsizes and 
prove that a linear convergence rate
is obtained. Numerical experiments for logistic regression problems are presented to
illustrate the impact of asynchrony upon convergence.

% Current literature propose to reduce the communication frequency by allowing the processors to run the descent algorithm locally and independent of each other and average the sequences regularly. This scheme has a major drawback: communication needs to be initiated periodically for every processor to perform averaging to ensure all processors have the same updated model. In this work, however, we allow the possibility of communication delays and some processors to have higher computation frequency than others. The advantage of this is first to reduce the synchronization penalty and second, to increase the implementation flexibility to handle the communication cost.  

\end{abstract}
\section{Introduction}

Asynchronous parallel optimization algorithms have gained attention in part due to increases in available data and use of parallel computation. 
These algorithms are used in large-scale machine learning  problems~\cite{lian2015asynchronous} and federated learning problems~\cite{bonawitz2019towards}. In control theory, similar applications arise in filtering  \cite{swerling71} and system identification \cite{aastrom1971system}, which lead to large optimization problems. % Similar applications also arise in power control~\cite{chang2014distributed},
% % robust sensor network control~\cite{rabbat2004distributed}
% signal processing~\cite{olshevsky2010efficient}, and robotics. 
Asynchronous algorithms are useful in parallel computing because they are not hindered by slow individual processors and they relax communication overhead compared to synchronized implementations.

This paper considers a class of optimization problems whose objective functions satisfy the Polyak-{\L}ojasiewicz (PL) condition, which is an inequality characterizing the curvature of some nonconvex functions~\cite{polyak1963gradient,karimi2016linear,yi2020primal}. The PL condition %can be regarded as a generalization of strong convexity in which 
requires neither convexity nor uniqueness of minimizers. Several important applications in machine learning have objective functions
that satisfy the PL condition; see~\cite{karimi2016linear} and references therein.
%, including least squares, rank-deficient least squares, logistic regression, and support vector machines;  
In control, PL functions can arise in state estimation and system identification \cite{swerling71, aastrom1971system},
which both use various forms of least squares. When such problems are rank deficient, they can fail to be strongly convex, though the PL condition still holds.
Thus, we expect the developments in this paper to be useful in state estimation and system ID, as well as other control problems which use optimization.

Recent work in~\cite{haddadpour2019local,yi2020primal} also studies this class of functions for a team of agents with local objective functions. That work uses an algorithmic model in which
each agent updates all decision variables and then agents average their iterates. Our algorithmic model is \emph{parallel}, in that each
decision variable is updated only by a single agent. 

For nonconvex problems, 
one way to accelerate classical gradient
descent algorithms is to use multiple processors to compute local gradients and update
their iterates using averages of gradients received from other processors.
For~$T$ iterations and~$n$ processors, 
this approach achieves $\mathcal{O} (\nicefrac{1}{\sqrt{nT}})$ convergence
for strongly convex functions and $\mathcal{O} (\nicefrac{1}{nT})$ for smooth nonconvex stochastic 
optimization \cite{stich2018local}. However, the proposed linear speedup can be difficult to attain 
in practice because of the communication overhead it incurs~\cite{stich2018local, haddadpour2019local}. 

We consider an alternative algorithmic model that tolerates longer delays under weaker assumptions.
The algorithm we consider is asynchronous parallel block coordinate descent (BCD).
Although such update laws
have been studied before~\cite{peng2016arock,bertsekas1989parallel}, this work is, to the best of our knowledge,
the first to connect it to the PL condition.

% [Kasra:Note that in comparison with classical data-parallel SGD,  talk about the different sequence argument here or no?]

\textbf{\textit{Contributions:}} The main results of this paper are:
\begin{itemize}
\item We show that the asynchronous block coordinate descent algorithm converges to a global
minimizer in linear time under the PL condition. 
Compared with
recent work~\cite{tang2018d,yu2019computation,haddadpour2019local}, 
we achieve the same convergence rate under more general assumptions on the cost function, network architecture, and communication requirements. To the best of our knowledge, this work is the first to establish a linear speed up with arbitrary (but bounded) delays in parallelized computations and communications when minimizing PL functions. 

\item We expand our results to show that the PL condition
is weaker than the so-called Regularity Condition (RC) that has seen wide use in the
data science community, e.g.,~\cite{candes2015phase}.
RC can be used to show that gradient descent converges to a minimizer 
at a linear rate~\cite{candes2015phase}, and functions that satisfy it have
been studied in~\cite{chi2019nonconvex}. We leverage this result to show that the
asynchronous block coordinate descent algorithm attains a linear convergence rate for
this class of functions as well. 
\item The asynchronous BCD algorithm is a standard algorithm, though the analysis and convergence results that we present for PL functions are entirely novel. 
Our work is closest to \cite{tseng1991rate}, which presents a block coordinate descent algorithm for objective functions that satisfy a form of the ``error-bound condition'' \cite{cannelli2020asynchronous}. We derive analogous convergence results, but under weaker assumptions and with a substantially simplified proof. In particular, we develop a novel proof strategy to leverage the PL property to prove convergence and derive a convergence rate tailored to PL functions.
\end{itemize}

This paper is organized as follows. Section \ref{sec:problemFormulation} provides a problem formulation. 
Section~\ref{sec:results} shows the linear convergence of our algorithm. 
We provide simulation results in Section~\ref{sec:casestudy}, and Section~\ref{sec:conclusions} concludes the paper.

\section{Background and Asynchronous Algorithm}\label{sec:problemFormulation}
This section presents the asynchronous parallel implementation of
coordinate descent
and assumptions we use to derive convergence rates. 
We use the notation~$[n] := \{1, \ldots, n\}$. 
%In particular, we allow
%arbitrarily large communication delays and 
%arbitrarily large differences in computation speeds among different processors provided that these delays obey some uniform upperbound.
\subsection{Optimization Problem}
We consider~$n$ processors jointly solving 
\begin{equation}\label{eq: minimize}
\min_{x \in \mathbb{R}^{m}} f(x),
\end{equation}
where~$f:\mathbb{R}^m\rightarrow \mathbb{R}$ is a continuously differentiable function and
satisfies the Polyak-{\L}ojasiewicz inequality:
\begin{dfn}\label{dfn:PL}\emph{(Polyak-{\L}ojasiewicz (PL) Inequality)} A function
satisfies the PL inequality if, for some $\mu>0$,
\begin{equation}
\frac{1}{2}\|\nabla f(z)\|^{2}\geq\mu\left(f(z)-f^*\right)
\quad
\textnormal{ for all } 
\quad
z \in \mathbb{R}^m,
\end{equation}
where~$f^* = \min_{x \in \mathbb{R}^m} f(x)$.
We say such an~$f$ is~$\mu$-PL or has the~$\mu$-PL property. 	\hfill $\triangle$
\end{dfn}
A~$\mu$-PL function has a unique global minimum value,
denoted by~$f^*$, and
the PL condition implies that every stationary point is a global minimizer. The $\mu$-PL property is implied by $\mu$-strong convexity, but it allows
for multiple minima and does not require convexity of any kind. 
For example,~$f(x) = x^2 + 3\sin^2(x)$ is non-convex and
satisfies the PL
inequality with~$\mu=1/32$.
It has also been
shown to be satisfied by
problems in signal processing and machine learning, including phase retrieval
\cite{candes2015phase}, some neural networks~\cite{zhou2017characterization}, matrix sensing, and matrix
completion \cite{bhojanapalli2016global}. 
We assume the following about~$f$.
%which, apart from the~$\mu$-PL property,
%are common in descent algorithms~\cite{bertsekas1989parallel}. 

\begin{assumption} \label{asmp:basic-f-assmp}
\leavevmode
\begin{enumerate}
    \item $f$ is~$\mu$-PL for some~$\mu > 0$.
    \item The set $\mathcal{X}^{* } \!=\! \{x^* \in \mathbb{R}^{m} \mid \nabla f(x^*)=0\}$ is nonempty and finite. \label{as:f2}
    \item $\nabla f(x)$ is $L$-Lipschitz continuous. In particular, \label{as:f3}
    \begin{equation}
        f(y) \leq f(x) + \nabla f(x)^T(y-x)
        + \frac{L}{2}\|y - x\|^2.
    \end{equation}
\end{enumerate}
\end{assumption}

\subsection{Asynchronous Parallel Block Coordinate Descent}
We decompose~$x$ via~${x = (x_1, \ldots, x_n)^T}$, 
where~$x_i \in \mathbb{R}^{m_i}$ and~$m = \sum_{i=1}^{n} m_i$.
Below, processor~$i$ computes updates only for~$x_i$. 
Define~$\nabla_i f = \frac{\partial f}{\partial x_i}$. 
%We use~$T^i \subseteq \mathbb{N}$ to denote
%the times at which processor~$i$ updates~$x_i$. 

Each processor stores a local copy of the decision variable~$x$. Due to asynchrony, 
these can disagree. We denote processor~$i$'s decision variable at time~$t$
by~$x^i(t)$. Processor~$i$ computes updates to~$x^i_i$ but not~$x^i_j$ for~$j \neq i$.
Instead, processor~$j$ updates~$x^j_j$ locally and transmits updated values
to processor~$i$. Due to asynchrony these values are delayed, and, in particular,~$x^i_j(t)$
can contain an old value of~$x^j_j$. We define~$\tau^i_j(t)$ to be the
time at which processor~$j$ originally computed the value that processor~$i$
has stored as~$x^i_j(t)$.
That is,~$\tau^i_j$ satisfies~$x^i_j(t) = x^j_j\big(\tau^i_j(t)\big)$. 
Clearly~$\tau^i_j(t) \leq t$ and we have
%\begin{equation}\label{key}
$x^{i}(t)=\Big(x^1_{1}\big(\tau_{1}^{i}(t)\big), \ldots, x^i_{i}\left(t\right),\ldots, x^n_{n}\big(\tau_{n}^{i}(t)\big)\Big)$. 
%\end{equation}
Below, we will also analyze the ``true" state of the network, which we define as
\begin{equation} \label{eq:xtrue}
x(t) = \big(x^1_1(t), x^2_2(t), \ldots, x^n_n(t)\big).
\end{equation}

We define~$T^i \subseteq \mathbb{N}$ as the set of times at which
processor~$i$ updates~$x^i_i$; agent~$i$ does not
actually know (or need to know)~$T^i$ because it is merely a tool
used for analysis. 
For all~$i \in [n]$ and stepsize $\gamma > 0$,
processor~$i$ executes
\begin{align} \label{eq:gd_updata}
x^i_{i}(t+1)&=
\begin{cases}
x^i_{i}(t)-\gamma \nabla_{i} f\left(x^{i}(t)\right) & t \in T^{i} \\
x^i_{i}(t) & \text { otherwise }
\end{cases}.
\end{align}

Our usage of the sets $T^i$ and time instants $\tau^i_j (t)$ is similar to~\cite[Chapter 6]{bertsekas1989parallel}: they are introduced for our use analytically and enable expression of asynchrony in algorithms that lack a common clock, though we emphasize that they need not be known to agents.
We assume that communication and computation delays are bounded, which
has been called \emph{partial asynchrony} in the literature~\cite{bertsekas1989parallel}. Formally, we have: 

\begin{assumption}\label{asmp:partial_async}
	There exists a positive integer~$B$ such that
	\begin{enumerate}
		\item For every~$i \in [n]$ and~$t \in \mathbb{N}$, at least one of the elements of the 
		set~$\{t, t+1, \dots, t+B-1\}$ is in~$T^i$.
		\item There holds~$t-B < \tau^i_j (t) \leq t$, 
		for all~$i,j \in [n]$,~$j \neq i$, and all~$t \in T^i$.
% 		\item There holds $\tau^i_i (t) = t$ for all $ i $ and $ t \in T^i $.
	\end{enumerate}	
\end{assumption}

We summarize the algorithm as follows. 

\begin{algorithm}\label{alg:GD}
\SetAlgoLined
\SetKwInOut{Input}{Input}
\SetKwInOut{Initialize}{Initialize}

\Input{Choose a stepsize $\gamma > 0$} 

\Initialize{$\{x^i\}_{i=1}^{n}$ for $n$ processors}
 \For{$t= 0,1, \ldots$, T}{
  \For{$i \in [n]$} {
  \eIf{$t\in {T^i}$}{
   Update: $x^i_{i}(t+1) = x^i_{i}(t)-\gamma \nabla_{i} f\left(x^{i}(t)\right)$
   }{
   Do not Update: $x^i_{i}(t+1) = x^i_{i}(t)$
  }
  \For{$j \in [n]\backslash\{i\}$}{
%   \begin{equation*}
      
%   \end{equation*}
%     x^i_{j}(t+1)=\left\{\begin{array}{ll}
%     x^i_{i}(\tau _{i}^{j}(t+1)) & j \text { receives } x_{i}^{i} \text { at } t+1 \\
%     x^i_{j}(t+1) & \text { otherwise }
%     \end{array}\right
\eIf{\textnormal{processor}~$i$ \textnormal{receives}~$x^i_j$ \textnormal{at time}~$t+1$}{
     $x_j^i(t+1) = x_{j}^{j}(\tau_{j}^{i}(t+1))$
     }{
    $x_j^i(t+1) = x_{j}^{i}(t)$
     }
    }
%\begin{equation*}
%x_{j}^{i}(t+1)=\begin{cases}
%x_{i}^{i}(\tau_{i}^{j}(t+1)) & j\text{ receives }x_{i}^{i}\text{ at }t+1\\
%x_{j}^{i}(t) & \text{otherwise}
%\end{cases}
%\end{equation*}

  }
}
\caption{Asynchronous BCD}
\end{algorithm}

% Stepsize bounds are critical to ensuring convergence, which is
% the focus of the next section. 

\section{Convergence Analysis}\label{sec:results}
We first prove linear convergence of Algorithm~1
under the PL condition.
Then we show that the Regularity Condition (RC)~\cite{candes2015phase} implies
the PL inequality and provide convergence guarantees for RC functions as well. 

\subsection{Convergence Under the PL Inequality}
Define
\begin{equation}
s_{i}(t):=\begin{cases}
-\nabla_{i}(f(x^{i}(t)) & t\in T^{i}\\
0 & \text{otherwise}
\end{cases}
\end{equation}
and concatenate the terms in $s(t) := [s_{1}(t)^T,\ldots,s_{n}(t)^T]^T$. 
%Let $\gamma \in (0,1)$ and take~$\eta > 0$ be large enough such that
%\begin{equation}\label{eq:eta}
%    f(x(t))-f^{*}\le \eta, \quad  \gamma^2 \sum_{\tau=t-B}^{t-1}\left\Vert s(\tau)\right\Vert ^{2}\le \eta
%\end{equation}
%for all~$t$.}

\begin{thm} \label{thm:main}  
Let~$f$ satisfy Assumption~\ref{asmp:basic-f-assmp}
and let Assumption~\ref{asmp:partial_async} hold. 
% Let~$\eta > 0$ be large enough \footnote{Note that the constant $\eta$ can be arbitrarily large and therefore we do not require neither $f(x)$ or $\nabla f(x)$ to be bounded.} such that
% \begin{equation}
%     f(x(t))-f^{*}\le \eta, \quad  \gamma^2 \sum_{\tau=t-B}^{t-1}\left\Vert s(\tau)\right\Vert ^{2}\le \eta
% \end{equation}
% for all~$t$.
%  where 
% \begin{align*} 
% \gamma  \le&\min \Bigg\{\frac{1}{\mu},
% \frac{1}{2L}\frac{1}{LnB+B+1},
% \frac{1}{L\frac{B}{2}\left(3n+1\right)+L+\mu+1},\\&
% \frac{\mu}{A_{2} Ln}\frac{1}{\left(1+L\left(B(n+1)+2\right)\right)\left(nBL^{2}+L+B+1\right)}
% \Bigg\},
% \end{align*}\label{eq:min_of_gamma}
% and the positive constants $A_1, A_2$ are defined as $A_1:=\frac{L}{2}nB(1+2L(\frac{B}{2}(n+1)+1)),$ $A_2:=2\mu+A_1\left(\frac{L}{2}nB+4LnB+8\right)$.
% \begin{align}
%  &A_1:=\frac{L}{2}nB(1+2L(\frac{B}{2}(n+1)+1)), \\& A_2:=2\mu+A_1\left(\frac{L}{2}nB+4LnB+8\right). 
% \end{align}
There exists $\gamma_0 \in (0, 1)$ such that for all $\gamma \in (0, \gamma_0)$ and~$x(t)$ as defined in~\eqref{eq:xtrue}, 
the sequence $\{x(t)\}_{t \in \mathbb{N}}$ generated by Algorithm~\ref{alg:GD} satisfies
\begin{equation}
f(kB) -f^{*}\le\left(1-\gamma\mu\right)^{k-1}\eta,\label{eq:induction_a}
\end{equation}
\begin{equation}
\gamma ^ 2 \sum_{\tau=(k-1)B}^{kB-1}\left\Vert s(\tau)\right\Vert ^{2}\le\left(1-\gamma\mu\right)^ {k-1}\eta \label{eq:induction_b}
\end{equation}
for some finite constant~$\eta > 0$ and for all $k=0,1,2,\dots$. 
\end{thm}
\noindent \emph{Proof:}
See Appendix. \hfill $\blacksquare$

This result generalizes standard results for strongly convex functions to the
case of~$\mu$-PL functions minimized under asynchrony.
If we consider centralized gradient descent for a~$\tau$-strongly convex
function, then~$B=1$ and~$\mu=\tau$, and we recover the classic linear rate
for strongly convex functions. 

\subsection{Convergence Under RC}
We next extend Theorem~\ref{thm:main} to 
objective functions that satisfy the Regularity Condition introduced in~\cite{candes2015phase}:

\begin{dfn} \emph{(Regularity Condition)} A function $f$ 
satisfies the
Regularity Condition $\text{RC}(\alpha,\beta)$ with~$\alpha, \beta > 0$, if
\begin{equation}
    \left\langle \nabla f(z),z-x^{*}\right\rangle \geq\frac{1}{\alpha}\|\nabla f(z)\|^{2}+\frac{1}{\beta}\left\Vert z-x^{*}\right\Vert ^{2}
\end{equation}
for all $z$, where $x^{*}$ is a  minimizer of~$f$. \hfill $\triangle$
\end{dfn}
We say such an~$f$ is~$\text{RC}(\alpha, \beta)$. 
This condition has appeared in machine learning and signal processing applications including  phase retrieval and matrix sensing~\cite{candes2015phase}. While it is simple to show that centralized gradient descent converges linearly for such functions, the convergence of parallel optimization algorithms under RC has received less attention, and we therefore extend our results to this case here. 
%We deepen our under    standing of the relationship between RC and PL inequality, and through a short extension of our results, we show that Algorithm~\ref{alg:GD} converges linearly to  minimizer of functions under RC.

Though elementary, we were unable to find the following lemma in the literature.

\begin{lemma}
Let~$f$ have a Lipschitz continuous gradient with Lipschitz constant~$L$. 
%For a function ~$f$ with Lipschitz-continuous gradient, there holds 
If~$f$ is~$\text{RC}(\alpha, \beta)$, then it is~$\frac{1}{\beta^2L}$-PL. 
%\begin{equation}
%    \textnormal{RC}(\alpha, \beta)\rightarrow \frac{1}{\alpha^2L}\textnormal{-PL}.
%\end{equation}
\end{lemma}
\noindent \emph{Proof:} 
Applying the Cauchy-Schwarz inequality to the RC definition, we write 
\begin{equation}
    \left\Vert \nabla f(z)\right\Vert \left\Vert z-x^{*}\right\Vert \geq\frac{1}{\alpha}\|\nabla f(z)\|^{2} + \frac{1}{\beta}\left\Vert z-x^{*}\right\Vert^{2}. 
\end{equation}
This gives  $\left\Vert \nabla f(z)\right\Vert \geq\frac{1}{\beta}\left\Vert z-x^{*}\right\Vert$. 
% and $\left\Vert z-x^{*}\right\Vert \geq\frac{1}{\alpha}\|\nabla f(z)\|$. 
%Thus,~$f$ satisfying~$\textnormal{RC}(\alpha,\beta)$ gives
%\begin{equation}\label{eq:alternativeRC}
%\frac{1}{\beta}\left\Vert z-x^{*}\right\Vert \le\left\Vert \nabla f(z)\right\Vert \le\alpha\left\Vert z-x^{*}\right\Vert .
%\end{equation}
Using Lipschitz continuity of the gradient (cf. Assumption~\ref{asmp:basic-f-assmp}.\ref{as:f3}),
for any~${x^* \in \mathcal{X}^*}$ and~$z \in \R^m$
we obtain
\begin{equation}
f(z)-f\left(x^{*}\right)\leq\frac{L}{2}\left\Vert z-x^{*}\right\Vert ^{2}\leq\frac{\beta^2 L}{2}\|\nabla f(z)\|^{2},
\end{equation}
and~$f$ satisfies the PL inequality with $\mu = \frac{1}{\beta^2L}$.  
$\hfill\blacksquare$

This lets us use Algorithm~\ref{alg:GD} for RC functions. 
\begin{thm}
Let~$f$ be~$\text{RC}(\alpha, \beta)$ and have Lipschitz gradient with constant~$L$. Let~$\mathcal{X}^*$ be finite, and let Assumption~\ref{asmp:partial_async} hold. 
Let $x^* \in \mathcal{X}^*$ be a minimizer of~$f$. 
For $\gamma \in (0,\gamma_0)$ as in Theorem~\ref{thm:main}, the sequence $\{x(t)\}_{t \in \mathbb{N}}$ generated by Algorithm~\ref{alg:GD} converges linearly to a fixed point $x^*$ as in Theorem~\ref{thm:main} with $\mu = \frac{1}{\beta^2L}$.
\end{thm}
\noindent \emph{Proof:} 
Immediately follows from Theorem~\ref{thm:main}. \hfill $\blacksquare$

\section{Case Study}\label{sec:casestudy}
%\textbf{Distributed Regularized Logistic Regression} 
We solve an $\ell_2$-regularized logistic regression problem using Algorithm~\ref{alg:GD}. This problem is strongly convex and differentiable, 
and therefore it satisfies the PL condition  \cite[Chapter 4]{aggarwal2020linear}.  We denote training feature vectors by $z^{(i)}\in \mathbb{R}^m$,
and we use $y^{(i)}\in \{0,1\}$ to denote their corresponding labels. The logistic regression objective function for $N$ observations is 
\begin{align}
    E(x)&=-\frac{1}{N}\left[\sum_{i=1}^{N}y^{(i)}\log(h_{i}(z^{(i)}))\right. \\
    &+(1-y^{(i)})\log(1-h_{i}(z^{(i)}))\Bigg]+\frac{\lambda}{2N}\left\Vert x\right\Vert ^{2}, 
\end{align}
where $h_{i}(x) = \frac{1}{1+e^{-x^{T}z^{(i)}}}$ is a sigmoid hypothesis function.

We conduct experiments on the Epsilon dataset using the above logistic regression model. 
The Epsilon dataset is a popular benchmark for large scale binary classification~\cite{haddadpour2019local},
and it consists of $400,000$ training samples and $100,000$ test samples. Each sample has a feature dimension of $m = 2000$. All data is preprocessed to mean zero, unit variance, and normalized to a unit vector. 

\begin{figure}
\begin{center}
\includegraphics[scale=0.068]{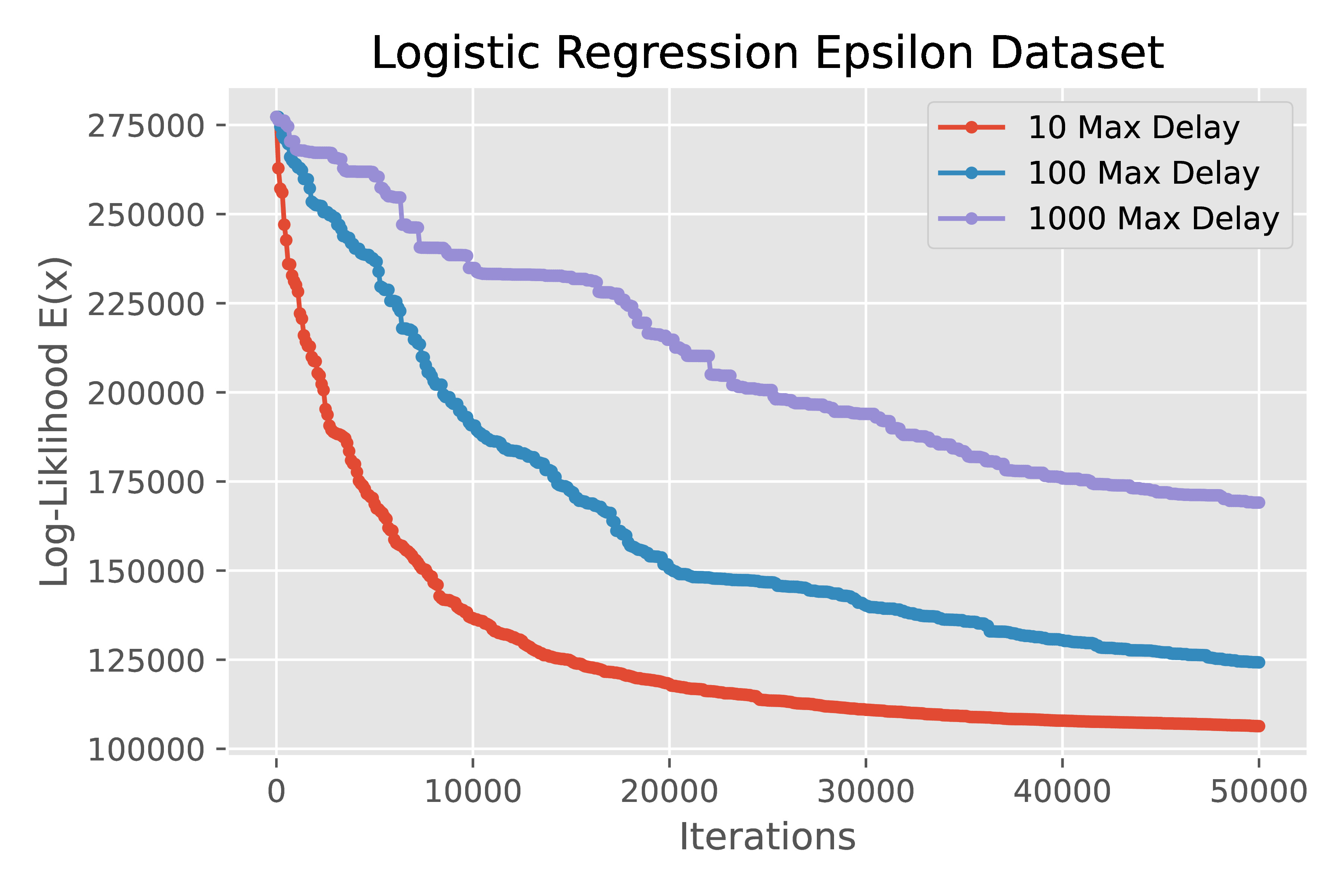}
\caption{Comparison between the convergence rate of BCD in Algorithm~\ref{alg:GD} with the maximum delay $B\in\{10,100, 1000\}$.}
\label{fig:costsWithDelay}
\end{center}
\end{figure}

We ran Algorithm~\ref{alg:GD} with $20$ processors with~$\gamma = 10^{-3}$
% \mh{Wait a minute. We need to make sure this satisfies all
% of our bounds in Theorem~1. If we're changing~$B$ here, shouldn't the stepsize also vary?}
and $\lambda = {10^{-2}}$. 
In three separate experiments, 
the communication delay for each processor is randomly generated
and bounded by $B=10$,~$B=100$, and~$B=1000$, respectively. 
The results of the experiment are shown in Figure~\ref{fig:costsWithDelay}. 

These results show
that Algorithm~\ref{alg:GD} converges linearly.
Indeed, we can observe that it converges with a slower rate as the communication delays increase, which
reflects the ``delayed linear'' nature of Theorem~\ref{thm:main}, which contracts toward
a minimizer by a factor of~$1-\gamma\mu$ every~$B$ timesteps. 

\section{Conclusions} \label{sec:conclusions}
We derived convergence rates of asynchronous coordinate decent parallelized among $n$ processors for 
functions satisfying the Polyak-{\L}ojasiewicz (PL) condition and those satisfying the Regularity Condition (RC). 
Future work includes
deriving similar convergence rates for stochastic settings.

%%%%%%%%%%%%%%%%%%%%%%%%%%%%%%%%%%%%%%%%%%%%%%%%%%%%%%%%%%%%%%%%%%%%%%%%%%%%%%%%

\bibliographystyle{IEEEtran}{}
\bibliography{sources}
%\newpage
%\clearpage
\section{Appendix}\label{sec:privacyImplementation}

We  begin with the following basic lemmas. 
\begin{lemma}\label{lem:partial_async}
For all $t\ge0$ and all~$i$, we have 
${
\left\Vert x^{i}(t)-x(t)\right\Vert \le\gamma\sum_{\tau=t-B}^{t-1}\left\Vert s(\tau)\right\Vert.  
}$
\end{lemma}
%\begin{IEEEproof}
\noindent \emph{Proof:} 
See Equation~(5.9) in~\cite[Section 7.5]{bertsekas1989parallel}. \hfill $\blacksquare$
%This follows from \cite[Section 7.5]{bertsekas1989parallel}. By the partial asynchrony Assumption \ref{asmp:partial_async}, we have
%\begin{align*}
%    \left\Vert x^{i}_j(t)-x_j(t)\right\Vert &=\left\Vert x_j\left(\tau_{j}^{i}(t)\right)-x_j(t)\right\Vert \\&=\gamma\left\Vert \sum_{\tau=\tau_{j}^{i}(t)}^{t-1}s_j(\tau)\right\Vert \leq\gamma\sum_{\tau=t-B}^{t-1}\left\Vert s_j(\tau)\right\Vert.
%\end{align*}

%This holds for all $j$ and the  lemma follows. $\hfill\blacksquare$
%\end{IEEEproof}

Next, %by using the partial asynchrony assumption and the Lipschitzian property of the gradient, 
we quantify the $B$-step decrease in the function value in the true state $x(t)$.
\begin{lemma}
\label{lem:B_step_difference}For all $t\ge0$, we have
\begin{align}%\label{eq:f(x(t+B))&-f(x(t))}
f&(x(t+B))-f(x(t))  \le\frac{L}{2}\gamma^{2}nB\sum_{\tau=t-B}^{t-1}\left\Vert s\left(\tau\right)\right\Vert ^{2}\nonumber \\
 & +\left(\gamma^{2}L\left(\frac{B}{2}(n+1)+1\right)-\gamma\right)\sum_{\tau=t}^{t+B-1}\|s(\tau)\|^{2}.\label{eq:f(x(t+B))-f(x(t))}
\end{align}
\end{lemma}
\noindent \emph{Proof:} 
From the Descent Lemma \cite[Section 3.2]{bertsekas1989parallel}, we have
\begin{align*}
&f(x(t+1))  =f(x(t)+\gamma s(t))\\
 & \le f(x(t))+\gamma\sum_{i=1}^{n}s_i(t)^T \nabla_{i}f(x(t))+L\gamma^{2}\left\Vert s(t)\right\Vert ^{2}.
\end{align*}
% Adding and subtracting~$\nabla_{i}f(x^{i}(t))$, we find
Adding~$( s_i(t)-s_i(t) )$ to $\nabla_{i}f(x(t))$
and applying the Lipschitz property of the gradient gives
\begin{align}
    &f(x(t+1))-f(x(t))  \\&
    \le L\gamma\sum_{i=1}^{n}\left\Vert \nabla_{i}f(x^{i}(t))\right\Vert \left\Vert x(t)-x^{i}(t)\right\Vert +(L\gamma^{2}-\gamma)\left\Vert s(t)\right\Vert ^{2}.
\end{align}
Employing Lemma~\ref{lem:partial_async}, and applying the inequality~$ab \leq \frac{1}{2}(a^2 + b^2)$ inside the sum, we find
% \begin{align*}
% &f(x(t+1))-f(x(t))  \\& \le-\gamma\sum_{i=1}^{n}\nabla_{i}f(x^{i}(t))^T\left[\nabla_{i}f(x(t))-\nabla_{i}f(x^{i}(t))\right]\\
%  & -\gamma\sum_{i=1}^{n}\nabla_{i}f(x^{i}(t))^T\nabla_{i}f(x^{i}(t))+L\gamma^{2}\left\Vert s(t)\right\Vert ^{2}. 
%  \end{align*}
%  By applying the Lipschitz property of the gradient and Lemma~\ref{lem:partial_async}, we obtain
%  \begin{align*}
%  &f(x(t+1))-f(x(t)) \\
 %\!\!\! \le L\gamma\sum_{i=1}^{n}\left\Vert \nabla_{i}f(x^{i}(t))\right\Vert \left\Vert x(t)-x^{i}(t)\right\Vert +(L\gamma^{2}-\gamma)\left\Vert s(t)\right\Vert ^{2}\\
 %&\!\! \!\le L\gamma\sum_{i=1}^{n}\left\Vert \nabla_{i}f(x^{i}(t))\right\Vert \gamma\sum_{\tau=t-B}^{t-1}\left\Vert s\left(\tau\right)\right\Vert +(L\gamma^{2}-\gamma)\left\Vert s(t)\right\Vert ^{2}\\
%  &\leq L\gamma^{2}\sum_{i=1}^{n}\sum_{\tau=t-B}^{t-1}\left\Vert \nabla_{i}f(x^{i}(t))\right\Vert \left\Vert s\left(\tau\right)\right\Vert +(L\gamma^{2}-\gamma)\left\Vert s(t)\right\Vert ^{2}.
%  \end{align*}
%  Then, applying the inequality~$ab \leq \frac{1}{2}(a^2 + b^2)$ inside the double sum, we find
 \begin{align}\label{eq:one_step_decrease}
 &f(x(t+1))-f(x(t)) \\
 %& \le\!\frac{L}{2}\!\gamma^{2}\!\sum_{i=1}^{n}\!\!\sum_{\tau=t-B}^{t-1}\!\!\!\!\left[\left\Vert \nabla_{i}f(x^{i}(t))\right\Vert ^{2}\!\!\!+\left\Vert s\left(\tau\right)\right\Vert ^{2}\right]\!\!+\!(L\gamma^{2}-\gamma)\!\left\Vert s(t)\right\Vert ^{2}\\
 &\leq \frac{L}{2}\!\gamma^{2}\left[B\left\Vert s\left(t\right)\right\Vert ^{2}+\!n\!\!\!\sum_{\tau=t-B}^{t-1}\left\Vert s\left(\tau\right)\right\Vert ^{2}\right]+(L\gamma^{2}-\gamma)\left\Vert s(t)\right\Vert ^{2}. 
\end{align}
The proof is completed by applying this
inequality successively to $t,t+1,\dots,t+B-1$, and summing them up. % we obtain
%\begin{align*}
%f(x(t+B))-&f(x(t))  \le\frac{L}{2}\gamma^{2}nB\sum_{\tau=t-B}^{t+B-1}\left\Vert s\left(\tau\right)\right\Vert ^{2}\\
% & +\left(\gamma^{2}L\left(\frac{B}{2}+1\right)-\gamma\right)\sum_{\tau=t}^{t+B-1}\|s(\tau)\|^{2},
%\end{align*}
%as desired. 
\hfill $\blacksquare$

The proof scheme up to this point closely follows \cite{bertsekas1989parallel,tseng1991rate}, which were not focused on PL functions. From this point on, we leverage 
the $\mu$-PL property of the objective function, and the following lemma and the rest of the results are new. 
The next lemma bounds the one-step change in~$x(t)$.
\begin{lemma}
\label{lem:square_difference}For all $t\ge0$, there holds
\begin{align*}
\left\Vert x(t+1)-x(t)\right\Vert ^{2} & \le\left(n^2B\gamma^{4}L^{2}+\gamma^{2}Ln\right)\sum_{\tau=t-B}^{t-1}\left\Vert s\left(\tau\right)\right\Vert ^{2}\\
 & +\left(\gamma^{2}+\gamma^{4}LnB\right)\|\nabla f(x(t))\|^{2}.
\end{align*}
\end{lemma}
\noindent \emph{Proof:} 
% Recall that
% \begin{equation}
% x(t) = \big(x^1_1(t)^T, \ldots, x^n_n(t)^T\big)^T.
% \end{equation}
% We also define
% \begin{equation}
% g_i(t) = \begin{cases}
% -\nabla_i f\big(x^i(t)\big) & \textnormal{agent~$i$ updates~$x^i_i$ at time~$t$} \\
% 0 & \textnormal{otherwise}
% \end{cases}
% \end{equation}
% and
% \begin{equation}
% s(t) = \big(g_1(t)^T, \ldots, g_n(t)^T\big)^T.
% \end{equation}
% Thus, using Assumption~\ref{asmp:partial_async}, we have
% \begin{equation} \label{eq:as2consequence}
% \left\Vert x^{i}(t)-x(t)\right\Vert \le\gamma\left\Vert \sum_{\tau=t-B}^{t-1}\left|g(\tau)\right|\right\Vert \leq \gamma\sum_{\tau=t-B}^{t-1}\left\Vert g(\tau)\right\Vert .
% \end{equation}
With~$x(t+1) = x(t) + \gamma s(t)$, we add and subtract~$\gamma \nabla f\big(x(t)\big)$ and apply the triangle inequality to find
%Using the triangle inequality, we write
\begin{align*}
\|x(t&+1)\!-\!x(t)\| \!
% =\! \|x(t\!+\!1)\!-\!x(t)\!+\!\gamma\nabla f(x(t)) \!-\! \gamma\nabla f(x(t))\| \\
%  & \le\|x(t\!+\!1)\!-\!x(t)\!+\!\gamma\nabla f(x(t))\| \!+\! \|\gamma\nabla f(x(t))\|\\
 & \le \|\gamma s(t)+\gamma\nabla f(x(t))\|+\|\gamma\nabla f(x(t))\|. 
\end{align*}
Squaring both sides, we expand to find
%\begin{multline}
%\left\Vert x(t+1)-x(t)\right\Vert^{2}  \le\|\gamma s(t)+\gamma\nabla f(x(t))\|^{2} \\ +\|\gamma\nabla f(x(t))\|^{2} +2\|\gamma s(t)+\gamma\nabla f(x(t))\|\|\gamma\nabla f(x(t))\|.
% \end{multline}
%Expanding the gradient terms into blocks gives
\begin{align}
 \|x(t &+ 1) - x(t)\|^2 \le \|\gamma\nabla f(x(t))\|^{2}  \\
 &+2\gamma \sum_{i=1}^{n}\|\gamma\nabla f(x(t))\|\left\Vert s_i(t) \!+\! \nabla_{i}f(x(t))\right\Vert \\
 & + \gamma^2\sum_{i=1}^{n} \|\nabla_{i}f(x(t)) + s_i(t)\|^{2}. 
 \end{align}
 Using the Lipschitz property of the gradient gives
 \begin{align}
\|x(t\!+\!1) \!&-\! x(t)\|^{2} \!\leq\! \|\gamma\nabla f(x(t))\|^{2} \!+\! \gamma^2L^2\sum_{i=1}^{n}\left\Vert x^{i}(t) \!-\! x(t)\right\Vert^{2} \\
 & +2\gamma L\sum_{i=1}^{n}\|\gamma\nabla f(x(t))\| \|x^{i}(t) - x(t)\|.
\end{align}
 Applying Lemma~\ref{lem:partial_async}, we expand to find
 \begin{align}
 &\|x(t+1) - x(t)\|^2 \!\leq
%  \gamma^{2}L^{2}\sum_{i=1}^{n}\left(\gamma\sum_{\tau=t-B}^{t-1}\left\Vert s(\tau)\right\Vert \right)^{2} \\
%  & +\|\gamma\nabla f(x(t))\|^{2} +2\gamma^{2}L\sum_{i=1}^{n}\sum_{\tau=t-B}^{t-1} \|\gamma\nabla f(x(t))\|\left\Vert s(\tau)\right\Vert 
n^2B\gamma^{4}L^{2}\!\!\!\sum_{\tau=t-B}^{t-1}\!\!\!\left\Vert s(\tau)\right\Vert ^{2}\!+\!\|\gamma\nabla f(x(t))\|^{2}\\
& +\gamma^{2}Ln\sum_{\tau=t-B}^{t-1}\left\Vert s(\tau)\right\Vert ^{2}+\gamma^{2}LnB\|\gamma\nabla f(x(t))\|^{2}, 
\end{align}
where we use~$(\sum_{i=1}^{K} y_i)^2 \leq K\sum_{i=1}^{K} y_i^2$
and~${ab \leq \frac{1}{2}(a^2 + b^2)}$. 
Rearranging completes the proof.  $\hfill\blacksquare$

The next lemma bounds the distance to minima. 
\begin{lemma}
\label{lem:B_step_distance_to_min} Take ${\gamma < \min \{\frac{2}{L}\frac{1}{Ln^2B+B+1}, \frac{1}{L\frac{B}{2}(n+1)+L + \frac{L}{2}nB} \}}$. For all $t\ge0$, we have
\begin{multline}\label{eq:B_step_distance_to_min}
  f(x(t+B))-f^*\le\left(1-C_{1}\right)\left(f(x(t))-f^*\right)\\+ (C_{2}+C_{3})\sum_{\tau=t-B}^{t-1}\left\Vert s\left(\tau\right)\right\Vert ^{2},
\end{multline}
where $C_1$, $C_2$, and~$C_3$ are positive constants defined as 
${
C_{1}=-\mu\left(\gamma^{4}L^{2}nB+\gamma^{2}\left(LB+L\right)-2\gamma\right)
}$, 
$
{C_{2}=\frac{1}{2}n^2B\gamma^{4}L^{3}+\frac{1}{2}\gamma^{2}L(L+1)n},
$ and 
$
{C_{3}:=\frac{L}{2}\gamma^{2}nB}
$.
\end{lemma}
\noindent \emph{Proof:} 
By Lipschitz continuity (cf. Assumption~\ref{asmp:basic-f-assmp}) we write 
\begin{align*}
f(x(t+1)) &- f(x(t)) \le \langle \nabla f(x(t)), \gamma s(t)\rangle \\
&+ \frac{L}{2}\left\Vert x(t + 1) - x(t)\right\Vert ^{2}\\
%  & =-\gamma\sum_{i=1}^{n}\nabla_{i}f(x^{i}(t))^T\nabla_{i}f(x(t))+\frac{L}{2}\left\Vert x(t+1)-x(t)\right\Vert ^{2}\\
%  & =\gamma\sum_{i=1}^{n}\left(\nabla_{i}f(x(t))-\nabla_{i}f(x^{i}(t))\right)\nabla_{i}f(x(t))\\
%  & -\gamma\sum_{i=1}^{n}\left(\nabla_{i}f(x(t))\right)^{2}+\frac{L}{2}\left\Vert x(t+1)-x(t)\right\Vert ^{2}\\
 &\le \gamma L\sum_{i=1}^{n}\left\Vert x(t)-x^{i}(t)\right\Vert \left\Vert\nabla_{i}f(x(t))\right\Vert\\
 & -\gamma \left\Vert \nabla f(x(t))\right\Vert ^{2} +\frac{L}{2}\left\Vert x(t+1)-x(t)\right\Vert ^{2},
\end{align*}
% where in the third step we add an subtract $\nabla_{i}f(x(t))$ and
% distributed terms, and on the fourth step we use the Lipschitz continuity
% of the gradient. 
where we added $\nabla f(x(t))- \nabla f(x(t))$ to $\gamma s(t)$ and used Cauchy-Schwarz and the Lipschitz continuity of the gradient. Next, applying Lemma~\ref{lem:partial_async} and then using $ab\le\frac{1}{2}(a^{2}+b^{2})$ on the first term of the RHS and simplifying we get
\begin{align*}
&f(x(t+1))-f(x(t)) \le \frac{L}{2}\left\Vert x(t+1)-x(t)\right\Vert ^{2} \\ 
% &\le\frac{1}{2}\gamma^{2}L\sum_{i=1}^{n}\sum_{\tau=t-B}^{t-1}\left(\left\Vert s(\tau)\right\Vert ^{2}+\left\Vert\nabla_{i}f(x(t))\right\Vert^{2}\right)\\
% &-\gamma\sum_{i=1}^{n}\left(\nabla_{i}f(x(t))\right)^{2}+\frac{L}{2}\left\Vert x(t+1)-x(t)\right\Vert ^{2}\\
& + \left(\frac{\gamma^{2}}{2}LB-\gamma\right)\left\Vert \nabla f(x(t))\right\Vert ^{2}+\frac{\gamma^{2}}{2}Ln\sum_{\tau=t-B}^{t-1}\left\Vert s(\tau)\right\Vert ^{2}.
\end{align*}
Applying Lemma \ref{lem:square_difference}, we obtain
\begin{align*}
&f(x(t+1))-f(x(t)) \\& \le\left(\frac{1}{2}nB\gamma^{4}L^{3}+\frac{1}{2}\gamma^{2}L(L+1)n\right)\sum_{\tau=t-B}^{t-1}\left\Vert s(\tau)\right\Vert ^{2} \\&+\left(\frac{1}{2}\gamma^{4}L^{2}n^2B+\gamma^{2}\left(\frac{1}{2}LB+\frac{L}{2}\right)-\gamma\right)\left\Vert \nabla f(x(t))\right\Vert ^{2}.
\end{align*}
The fact that $\gamma < 1 $ and the first term in the definition of $\gamma$ ensure that
$
{\frac{1}{2}\gamma^{4}L^{2}n^2B+\gamma^{2}\left(\frac{1}{2}LB+\frac{L}{2}\right)-\gamma\le0}
$.
Using this and the PL inequality, we have
% \begin{multline}
% f(x(t+1))-f(x(t))\le-C_{1}\left(f(x(t))-f\left(x^{*}\right)\right)\\+C_{2}\sum_{\tau=t-B}^{t-1}\left\Vert s(\tau)\right\Vert ^{2}
% \end{multline}
\begin{align}
f(x(t+1))-f^{*}&\le(1-C_{1})\left(f(x(t))-f^*\right)\\&+C_{2}\sum_{\tau=t-B}^{t-1}\left\Vert s(\tau)\right\Vert ^{2}, \label{eq:c3_c2_eq-1}
\end{align}
where we have also added $-f^{*}$ to both sides. For a small enough $\gamma$, repeating the argument in Lemma~\ref{lem:B_step_difference} from $t+1$ to $t+B$, we find
% \begin{multline}
% f(x(t+B))-f(x(t+1))\le\frac{L}{2}\gamma^{2}nB\sum_{\tau=t-B}^{t-1}\left\Vert s\left(\tau\right)\right\Vert ^{2}\\ + \left(\gamma^{2}L\left(\frac{B}{2}(n+1)+1\right)-\gamma\right)\sum_{\tau=t}^{t+B-1}\|s(\tau)\|^{2}.
% \end{multline}
% This allows us to write the equation
% above as 
\begin{equation}\label{eq:f(x(t+B))-f(x(t+1))}
    -f(x(t+1))\le-f(x(t+B))+\frac{L}{2}\gamma^{2}nB\sum_{\tau=t-B}^{t-1}\left\Vert s\left(\tau\right)\right\Vert ^{2}.
\end{equation}
Using~\eqref{eq:f(x(t+B))-f(x(t+1))}  in~\eqref{eq:c3_c2_eq-1}, completes the proof. \hfill $\blacksquare$

To streamline further analyses, we define 
${
    \alpha(t):=f(x(t))-f^{*},
}$
${
    \beta(t):=\gamma^2 \sum_{\tau=t-B}^{t-1}\left\Vert s(\tau)\right\Vert ^{2}.
}$
We next find upper bounds for $\beta (t)$ and $\beta(t-B)$. We substitute $t-B$ for $t$ and $t-2B$ for $t$ consecutively in~\eqref{eq:f(x(t+B))-f(x(t))}, and by rearranging the terms, we obtain 
\begin{equation}\label{eq:beta(t)bound}
   \gamma^{-2}\beta(t)\le\frac{1}{C_{4}}\left(\alpha(t-B)+\frac{L}{2}nB\beta(t-B)-\alpha(t)\right),
\end{equation}
\begin{multline}\label{eq:beta(t-B)bound}
    \gamma^{-2}\beta(t-B)\le\frac{1}{C_{4}}\bigg(\alpha(t-2B)\\+\frac{L}{2}nB\beta(t-2B)-\alpha(t-B)\bigg),
\end{multline}
where ${C_{4}:=\gamma-\gamma^{2}L\left(\frac{B}{2}(n+1)+1\right)} \in \mathbb{R}^+$. Next, we employ~\eqref{eq:beta(t)bound} to bound $\beta(t)$ in Lemma~\ref{lem:B_step_distance_to_min}, which gives
\begin{align}\label{eq:sub_magic}
\alpha(t+B)\le&\left(1-C_{1}-\frac{C_{2}+C_{3}}{C_{4}}\right)\alpha(t)+\\&\frac{C_{2}+C_{3}}{C_{4}}\left(\alpha(t-B)+\frac{L}{2}nB\beta(t-B)\right).
\end{align}
Using~\eqref{eq:beta(t-B)bound}, we bound the term $\beta(t-B)$ in~\eqref{eq:sub_magic}, and by  simplifying the terms, we obtain 
\begin{align}\label{eq:magic}
    &\alpha(t+B)\le\left(1-C_{1}-\frac{C_{2}+C_{3}}{C_{4}}\right)\alpha(t)
    + \left(\frac{C_{2}+C_{3}}{C_{4}}\right)\cdot \nonumber \\& 
    \biggl[(1-\frac{C_{3}}{C_{4}})\alpha(t-B)+\frac{C_{3}}{C_{4}}\alpha(t-2B)
    % \\&
    +\frac{L}{2}nB\frac{C_{3}}{C_{4}}\beta(t-2B)\biggr].
\end{align}
Next, by applying Lemma~\ref{lem:B_step_difference}, at ${t, t+B, t+2B,\dots}$ and summing up those equations, we obtain 
\begin{align}
    \alpha(t+kB) &\le  \alpha(t) +  \left(-C_{4}+C_{3}\right)\sum_{l=1}^{k-1}\sum_{\tau=t+\left(l-1\right)B}^{t+lB-1}\left\Vert s\left(\tau\right)\right\Vert ^{2}\\& -C_{4}\sum_{\tau=t+\left(k-1\right)B}^{t+kB-1}\|s(\tau)\|^{2}+C_{3}\sum_{\tau=t-B}^{t-1}\left\Vert s\left(\tau\right)\right\Vert ^{2}.
\end{align}
Using $\gamma < \frac{1}{L\frac{B}{2}(n+1)+L + \frac{L}{2}nB}$ makes  $\left(C_{3}-C_{4}\right)<0$, and therefore we get 
\begin{align}\label{eq:t+kb_before_lim}
&\alpha(t+kB)\le\alpha(t)+\frac{L}{2}nB\beta(t)\\&
\!\!\!+\!\left(L(\frac{B}{2}(n+1)+1)-\gamma^{-1}+\frac{L}{2}nB\right)\beta(t+B).
\end{align}
By the PL inequality, $f(x(t)) \ge f^*$ for all $x(t)$ and therefore ${\lim\inf_{t\rightarrow\infty}f(x(t))\ge f^*}$, which gives ${\alpha(t+kB) \ge 0 }$ for all $k$. From this fact and \eqref{eq:t+kb_before_lim}  we obtain 
\begin{align}\label{eq:magic2}
\beta(t+B)\le\frac{\alpha(t)+\frac{L}{2}nB\beta(t)}{\gamma^{-1}-L\left(\frac{B}{2}(n+1)+1\right)-\frac{L}{2}nB}.
\end{align}
Finally, using~\eqref{eq:magic} and~\eqref{eq:magic2} we are able to show the main result on asynchronous linear convergence of BCD.

\indent \emph{Proof of Theorem~1:} 
For sufficiently small~$\gamma$, each step of Algorithm~1 decreases the objective function~$f(x)$, e.g.,~\eqref{eq:one_step_decrease}. Moreover, since the solution set is finite (from Assumption~\ref{asmp:basic-f-assmp}) and $f(x)$ is bounded from below due to the PL inequality, the sequence $\{x(t)\}_{t \in \mathbb{N}}$ generated by Algorithm~1 converges. In particular, it is bounded.

We proceed to show the linear convergence by induction on $k$. Take the scalar $\eta$ large enough such that 
\begin{equation}
    \alpha(t),\; \alpha(t+B),\; \beta(t),\; \beta(t+B)\leq \eta;
\end{equation}
such an~$\eta$ exists because~$x(t)$ is bounded, and by continuity so are~$\alpha(t)$ and~$\beta(t)$. 
Then \eqref{eq:induction_a} and \eqref{eq:induction_b} hold for $k=0$ and $k=1$. We show that if $\alpha(t+kB)\le \left(1-\gamma\mu\right)^{k-1}\eta$ and $\beta(t+kB)\le \left(1-\gamma\mu\right)^{k-1}\eta$, then ${\alpha(t+(k+1)B)\le \left(1-\gamma\mu\right)^{k}\eta}$ and ${\beta(t+(k+1)B)\le \left(1-\gamma\mu\right)^{k}\eta}$ for $k \ge 1$.
By this induction hypothesis, 
\eqref{eq:magic} is written as 
\begin{align*}
    \alpha(&t+(k+1)B)\le(1-C_{1}-\frac{C_{2}+C_{3}}{C_{4}})\eta(1-\gamma \mu)^{k-1}\\&
    +\frac{C_{2}+C_{3}}{C_{4}}\Bigg((1-\frac{C_{3}}{C_{4}})\eta(1-\gamma \mu)^{k-2}+\frac{C_{3}}{C_{4}}\eta(1-\gamma \mu)^{k-3}\\&
    +(\frac{L}{2}nB)\frac{C_{3}}{C_{4}}\eta(1-\gamma \mu)^{k-3}\Bigg). 
\end{align*}
Factoring out the term $\left(1-\gamma \mu \right)^{k-1}\eta$,  using the inequality $\left(1-\gamma \mu \right)^{-1} <  (1+2\gamma \mu)$ and its square, and  simplifying the above equation we get
\begin{align*}
    \alpha&(t\!+\!(k\!+\!1)B) \le 
    \Biggl[1\!-\!C_{1} \!+\! \frac{C_{2} \!+\! C_{3}}{C_{4}} \Biggl(\!2\gamma \mu-\frac{C_{3}}{C_{4}}\left(1 \!+\! 2\gamma \mu\right)\\& +\frac{C_{3}}{C_{4}}(1+\frac{L}{2}nB)(1+4\gamma^{2}\mu^{2}+4\gamma \mu)\Biggr)\Biggr]\eta\left(1-\gamma \mu\right)^{k-1}. 
\end{align*}
Next, take $\gamma < \frac{1}{\mu}$ which gives $1+4\gamma \mu(1+\gamma \mu) < 1+8\gamma \mu$. We use that and the inequality, $\left(1-\gamma L\left(\frac{B}{2}(n+1)+1\right)\right)^{-1} < 1+2\gamma L\left(\frac{B}{2}(n+1)+1\right)$ and expand to get
\begin{align*}
&\alpha(t+(k+1)B) \le\Biggl[1+\mu\left(\gamma^{4}L^{2}nB+\gamma^{2}\left(LB+L\right)-2\gamma\right)\\&
+(1+2\gamma L(\frac{B}{2}(n+1)+1)) \, \cdot \\&
(\frac{L}{2}\gamma nB+\frac{1}{2}n^2B\gamma^{3}L^{3}+\frac{1}{2}\gamma L(L+1)n) \, \cdot\\&
(2\gamma \mu+\frac{C_{3}}{C_{4}}(\frac{L}{2}nB+8\gamma \mu+LnB4\gamma \mu))\Biggr]\eta\left(1-\gamma \mu\right)^{k-1}.
\end{align*}
Define the constants $A_1, A_2 \in \mathbb{R}^{+}$ by $A_1:=\frac{L}{2}nB(1+2L(\frac{B}{2}(n+1)+1)),$ $A_2:=2\mu+A_1\left(\frac{L}{2}nB+4LnB+8\right)$.
Note that from $\gamma < 1$ we have $\frac{C_{3}}{C_{4}} < \gamma A_1$. Simplifying gives 
\begin{align*}
\alpha(&t+(k+1)B)\le\Biggl[1-2\gamma\mu 
	+\mu\gamma^{2}\left(L^{2}nB+LB+L\right) \\&
	+\gamma^{2}\left(A_{2}+2A_{2}L(\frac{B}{2}(n+1)+1)\right) \cdot \\&
	\left(\frac{L}{2}nB+\frac{1}{2}n^2BL^{3}+\frac{1}{2}L(L+1)n\right)\Biggr]\eta\left(1-\gamma \mu\right)^{k-1}.
\end{align*}
Taking $\gamma$ according to
\begin{align*} 
\gamma  \le &\min \Bigg\{
\frac{1}{2L}\frac{1}{LnB+B+1},\\&
\frac{\mu}{A_{2} Ln}\frac{1}{\left(1+L\left(B(n+1)+2\right)\right)\left(nBL^{2}+L+B+1\right)}
\Bigg\}
\end{align*}
allows us to upper bound all of the above terms to get 
\begin{equation}
    \alpha(t+(k+1)B)\le\left(1-\gamma\mu\right)^{k}\eta. 
\end{equation}
This completes the first part of the proof. From the induction hypothesis, \eqref{eq:magic2} gives
\begin{equation}
   \beta(t+(k+1)B) \le \frac{\eta\left(1-\gamma \mu \right)^{k-1}+\frac{L}{2}nB\eta\left(1-\gamma \mu \right)^{k-1}}{\left(\gamma^{-1}-L\left(\frac{B}{2}(n+1)+1\right)-\frac{L}{2}nB\right)}. 
\end{equation}
Taking $\gamma$ according to  
$
  \gamma < \frac{1}{L\frac{B}{2}\left(3n+1\right)+L+\mu+1}
$
therefore gives
\begin{equation}
    \beta(t+(k+1)B)\le\left(1-\gamma\mu\right)^{k} \eta.
\end{equation}
This holds for all $t$. The value of $\gamma_0$ can be found by combining the upper bounds on admissible stepsizes $\gamma$ in Lemma~\ref{lem:B_step_distance_to_min} and those used above. $\hfill \blacksquare$

%\newpage
%\red{
%\begin{enumerate}
%\item Verify that all inequality arguments are valid. For example, suppose~$x_1 \leq y \leq x_2$, and suppose we want to derive a sufficient
%condition to have~$y \leq z$. It is sufficient to have~$x_2 \leq z$, but it is not sufficient to have~$x_1 \leq z$. 
%\item Make sure every symbol is defined, ideally before it is used. Pay special attention to network-level terms (like~$H$,~$V$, etc.). 
%\item Look at all margin comments and make the appropriate changes.
%\item Only use~$\tilde{y}$ to denote the privatized output. Specifically, only write~$\Delta_{\ell_2} y$ and don't write~$\Delta_{\ell_2}\tilde{y}$
%\item Put all assumptions in the theorem statement, not in the proof
%\end{enumerate}
%}

\end{document}